\documentclass[11pt]{article}
\usepackage{latexsym,amsfonts,amsmath,graphics}
\usepackage{epsfig}\usepackage{mathrsfs}

\DeclareMathOperator{\betafcn}{B}
\DeclareMathOperator{\betafcnReg}{I}
\DeclareMathOperator{\dist}{d}

\DeclareMathOperator{\dd}{d}
\newtheorem{theorem}{Theorem}

\newtheorem{proposition}{Proposition}

\newenvironment{proof}{\begin{trivlist}
    \item[\hskip\labelsep{\it Proof.}]}{$\hfill\Box$\end{trivlist}}

\newcommand{\bsa}{\boldsymbol{a}}

\newcommand{\bsx}{\boldsymbol{x}}

\newcommand{\bsp}{\boldsymbol{p}}
\newcommand{\bst}{\boldsymbol{t}}
\newcommand{\bsz}{\boldsymbol{z}}

\newcommand{\bsy}{\boldsymbol{y}}

\allowdisplaybreaks

\begin{document}

\title{\scshape A simple Proof of Stolarsky's \\ Invariance Principle\footnote{MSC: primary 41A30; secundary 11K38, 41A55.}}

\author{Johann S. Brauchart\thanks{\noindent The author is supported by an {\sc APART}-Fellowship of the Austrian Academy of Sciences.} and Josef Dick\thanks{The author is supported by an Australian Research Council Queen Elizabeth 2 Fellowship.} \\ School of Mathematics and Statistics,\\
University of New South Wales, \\
Sydney, NSW, 2052, Australia \\ \texttt{j.brauchart@unsw.edu.au}, \\ \texttt{josef.dick@unsw.edu.au}}

\date{}
\maketitle

\begin{abstract}
Stolarsky 
[Proc. Amer. Math. Soc. 41 (1973), 575--582] 
showed a beautiful relation that balances the sums of distances of points on the unit sphere and their spherical cap $\mathbb{L}_2$-discrepancy to give the distance integral of the uniform measure on the sphere a potential-theoretical quantity (Bj{\"o}rck [Ark. Mat. 3 (1956), 255--269]). Read differently it expresses the worst-case numerical integration error for functions from the unit ball in a certain Hilbert space setting in terms of the $\mathbb{L}_2$-discrepancy and vice versa (first author and Womersley [Preprint]). 
In this note we give a simple proof of the invariance principle using reproducing kernel Hilbert spaces. 
\end{abstract}

\section{Introduction}

We consider the unit sphere 
\begin{equation*}
\mathbb{S}^d = \left\{ \boldsymbol{z} =
(z_1,\ldots, z_{d+1}) \in \mathbb{R}^{d+1}: \|\boldsymbol{z}\| =
\sqrt{z_1^2+\cdots + z_{d+1}^2} = 1 \right\}
\end{equation*}
embedded in the Euclidean space $\mathbb{R}^{d+1}$, $d \geq 2$.
Let $f:\mathbb{S}^d \to\mathbb{C}$ be a continuous function. 
Then we approximate the integral $\int_{\mathbb{S}^d} f(\bsx) \, \dd  \sigma_d(\bsx)$, 
where $\sigma_d$ is the normalized Lebesgue surface area measure on $\mathbb{S}^d$ ($\int_{\mathbb{S}^d} \dd  \sigma_d = 1$), by an equal weight numerical integration rule 
\begin{equation} \label{eq:Q.N}
Q_N(f) = \frac{1}{N} \sum_{k=0}^{N-1} f(\bsz_k)
\end{equation}
where $\bsz_0,\ldots, \bsz_{N-1} \in \mathbb{S}^d$ are the integration nodes on the sphere. In order to analyze the integration error committed by the approximation, we define a worst-case error by
\begin{equation*}
e(\mathcal{H},Q_N) = \sup_{f \in \mathcal{H}, \|f\| \le 1} \left|\int_{\mathbb{S}^d} f(\bsx) \,\mathrm{d} \sigma_d(\bsx) - \frac{1}{N} \sum_{k=0}^{N-1} f(\bsz_k) \right|,
\end{equation*}
where $\mathcal{H}$ denotes a normed function space with norm $\|\cdot \|$. The rate of decay of the worst-case error depends on the function space and the integration nodes. For a fixed function space, the worst-case error can serve as a quality criterion for different sets of integration nodes, meaning that the performance of a set of $N$ quadrature points $\bsz_0,\ldots, \bsz_{N-1}$ can be compared to another set of $N$ quadrature points by comparing the corresponding worst-case errors. Generally, this only means that the integration error is smaller, but sometimes the worst-case error allows also a geometrical interpretation.

Another quality criterion for points on the sphere exploits the potential energy, or more generally, the Riesz $s$-energy of configurations of points modeling unit charges which are thought to interact through a potential $1/\| \cdot \|^s$ ($s\neq0$), where $\|\cdot\|$ denotes the Euclidean distance.(We refer the reader to the survey papers \cite{HaSa2004} and \cite{SaKu1997} and for {\em universally optimal configurations} to \cite{CoKu2007}.) A particular instance is the (normalized) sum of distances ($s=-1$) 
\begin{equation*}
\frac{1}{N^2} \sum_{k,\ell =0}^{N-1} \|\bsz_k - \bsz_\ell \|, \qquad \bsz_1, \dots, \bsz_N \in \mathbb{S}^d.
\end{equation*}
It is well-known from potential theory (see Bj{\"o}rck \cite{Bj1956}) that this (discrete) sum of distances of optimal $N$-point configurations approaches the associated (continuous) distance integral of the uniform measure $\sigma_d$ on $\mathbb{S}^d$ as $N \to \infty$. In fact, any sequence of $N$-point systems with this property turns out to be 'asymptotically uniformly distributed'; that is, the discrete probability measure obtained by placing equal charges at the points tends to the uniform measure in the weak-star sense. The difference
\begin{equation}\label{eq_intdoublesum}
\int_{\mathbb{S}^d} \int_{\mathbb{S}^d} \|\bsz - \bsx\| \,\mathrm{d} \sigma_d(\bsz) \,\mathrm{d} \sigma_d(\bsx) - \frac{1}{N^2} \sum_{k,\ell=0}^{N-1} \|\bsz_k-\bsz_\ell \|
\end{equation}
measuring the deviation between theoretical and empirical $(-1)$-energy quantifies the quality of points on the sphere (and, indirectly, their uniform distribution) using energy.
%
It should be mentioned that the upper bound of correct order $N^{-1-1/d}$ for \eqref{eq_intdoublesum} (for optimal configurations) was obtained by Stolarksy~\cite{St1973} using his invariance principle and a result of Schmidt~\cite{Sch1969} on the discrepancy of spherical caps. The correct-order lower bound ($N^{-1-1/d}$) was established by Beck~\cite{Be1984} using his Fourier transform technique.

The spherical cap discrepancy measures the maximum deviation between theoretical and empirical distribution with respect to spherical caps as test sets. It can be used to compare point sets on the sphere with respect to their distribution properties. To introduce the concept of spherical cap discrepancy, we require some notation. A spherical cap centered at $\bsx \in \mathbb{S}^d$ with 'height' $t \in [-1,1]$ is the set
\begin{equation*}
C(\bsx;t) = \{\bsz \in \mathbb{S}^d: \langle \bsx,\bsz \rangle \ge
t\}.
\end{equation*}
The family of all spherical caps is denoted by 
\begin{equation*}
\mathscr{C} = \{C(\bsx;t): \bsx \in \mathbb{S}^d, t \in [-1,1]\}.
\end{equation*}

For a set $J \subseteq \mathbb{R}^{d+1}$ we define the indicator function
\begin{equation*}
1_{J}(\bsx) = 
\begin{cases}
1 & \text{if $\bsx \in J$,} \\
0 & \text{otherwise.}
\end{cases}
\end{equation*}
For a measurable set $J \subseteq \mathbb{S}^d$ let
\begin{equation*}
\sigma_d(J) = \int_{\mathbb{S}^d} 1_{J}(\bsx) \,\dd  \sigma_d(x).
\end{equation*}
Using spherical caps we can define another quality criterion for points on the sphere $\mathbb{S}^d$ in terms of their distribution properties. One such criterion is the spherical cap $L_2$-discrepancy, which is given by
\begin{equation*}
L_2(P) = \left(\int_{-1}^1 \int_{\mathbb{S}^d} \left|\sigma_d(C(\bsz;t)) - \frac{1}{N} \sum_{k=0}^{N-1} 1_{C(\bsz;t)}(\bsz_k) \right|^2 \,\mathrm{d} \sigma_d(\bsz) \,\mathrm{d} t \right)^{1/2},
\end{equation*}
where $P = \{\bsz_0,\ldots, \bsz_{N-1}\}$.

We considered three, seemingly different, measures which can be applied to point sets on the sphere. It turns out that in some instances, the three measures are related to each other. Stolarksy's insight~\cite{St1973} was that the sum of distances of points on the sphere and the spherical cap discrepancy coincide. On the other hand, also the sum of distances of points on the sphere and the worst-case error coincide for a certain choice of function space, see \cite{BrWo2010d_pre} and also Sloan and Womersley~\cite{SlWo2004} regarding a generalized discrepancy of Cui and Freeden~\cite{CuFr1997}. 

In this paper we give a simple proof of these results based on reproducing kernel Hilbert spaces. We also provide some generalizations which follow from our approach.

\section{Reproducing kernel Hilbert space}

We define a reproducing kernel Hilbert space using the general approach of \cite[Ch.~9.6]{NoWo2010}.

For $\bsx,\bsy \in \mathbb{S}^d$ we define the function
$K_{\mathscr{C}}: \mathbb{S}^d \times \mathbb{S}^d \to \mathbb{R}$
by
\begin{equation} \label{eq:reproducing.kernel}
K_{\mathscr{C}}(\bsx,\bsy) = \int_{-1}^1 \int_{\mathbb{S}^d} 1_{C(\bsz;t)}(\bsx) 1_{C(\bsz;t)}(\bsy) \,\dd \sigma_d(\bsz) \,\dd t.
\end{equation}
Since $1_{C(\bsz;t)}(\bsx) = 1_{C(\bsx;t)}(\bsz)$, we also have
\begin{equation*}
K_{\mathscr{C}}(\bsx,\bsy) = \int_{-1}^1 \int_{\mathbb{S}^d} 1_{C(\bsx;t)}(\bsz) 1_{C(\bsy;t)}(\bsz) \,\dd \sigma_d(\bsz) \,\dd t.
\end{equation*}

The function is obviously symmetric, i.e. we have
$K_{\mathscr{C}}(\bsx,\bsy) = K_{\mathscr{C}}(\bsy,\bsx)$. Further,
let $a_0,\ldots, a_{N-1} \in \mathbb{C}$ and $\bsx_0,\ldots,
\bsx_{N-1} \in \mathbb{S}^d$. Then we have
\begin{align*}
\sum_{k,\ell=0}^{N-1} a_k \overline{a_\ell} K_{\mathscr{C}}(\bsx_k,\bsx_\ell) 
&= \int_{-1}^1 \int_{\mathbb{S}^d} \sum_{k,\ell=0}^{N-1} a_k
\overline{a_\ell} 1_{C(\bsx;t)}(\bsx_k) 1_{C(\bsx;t)}(\bsx_\ell)
\,\dd  \sigma_d(\bsx) \,\dd  t \\
&= \int_{-1}^1 \int_{\mathbb{S}^d} \left| \sum_{k=0}^{N-1} a_k
1_{C(\bsx;t)}(\bsx_k) \right|^2 \,\dd  \sigma_d(\bsx)
\,\dd  t \\ &\geq 0.
\end{align*}

Thus, the function $K_{\mathscr{C}}$ is symmetric and positive
definite. By \cite{Ar1950}, this implies that $K_{\mathscr{C}}$ is a
reproducing kernel. It is also shown in \cite{Ar1950} that a reproducing kernel uniquely defines a Hilbert space of functions with a certain inner product. Let $\mathcal{H}_{\mathscr{C}} =
\mathcal{H}(K_{\mathscr{C}}, \mathbb{S}^d)$ denote the corresponding
reproducing kernel Hilbert space of functions
$f:\mathbb{S}^d\to\mathbb{R}$ with reproducing kernel
$K_{\mathscr{C}}$.

We consider now functions $f_1,f_2:\mathbb{S}^d\to\mathbb{C}$ which
permit a certain integral representation. Let $g_1, g_2:
\mathbb{S}^d \times [-1,1] \to\mathbb{C}$ with $g_1,g_2 \in
\mathbb{L}_2(\mathbb{S}^d \times [-1,1])$ and
\begin{equation} \label{f_rep}
f_i(\bsx) = \int_{-1}^1 \int_{\mathbb{S}^d} g_i(\bsz;t)
1_{C(\bsz;t)}(\bsx) \,\dd \sigma_d(\bsz) \,\dd  t, \qquad i = 1, 2.
\end{equation}
Notice that for any fixed $\bsy \in \mathbb{S}^d$
the function $K_{\mathscr{C}}(\cdot,\bsy)$ also is of this form,
where the function $g$ is given by $1_{C(\bsz;t)}(\bsy)$ (considered
as a function of $\bsz$ and $t$ and where $\bsy$ is fixed). For
functions of this form we can define an inner product by
\begin{equation} \label{inner_product}
\langle f_1,f_2 \rangle_{K_{\mathscr{C}}} = \int_{-1}^1
\int_{\mathbb{S}^d} g_1(\bsz;t) \overline{g_2(\bsz;t)} \,\mathrm{d}  
\sigma_d(\bsz) \,\mathrm{d}  t.
\end{equation}
Let $\bsy \in \mathbb{S}^d$ be fixed. With this definition we obtain
\begin{equation*}
\langle f_1, K_{\mathscr{C}}(\cdot,\bsy) \rangle_{K_\mathscr{C}} =
\int_{-1}^1 \int_{\mathbb{S}^d} g_1(\bsz;t) 1_{C(\bsz;t)}(\bsy)
\,\dd \sigma_d(\bsz) \,\dd  t = f_1(\bsy).
\end{equation*}
By \cite{Ar1950}, the inner product in $\mathcal{H}_{\mathscr{C}}$ is
unique. Therefore, functions $f_i$, which are given by \eqref{f_rep} and for which $\langle f_i, f_i \rangle_{K_\mathscr{C}} < \infty$,
are in $\mathcal{H}_{\mathscr{C}}$ and \eqref{inner_product} is an
inner product for those functions in $\mathcal{H}_{\mathscr{C}}$.

Consider now the reproducing kernel $K_{\mathscr{C}}$. We have
\begin{equation*}
\int_{-1}^1 1_{C(\bsz;t)}(\bsx) 1_{C(\bsz;t)}(\bsy) \,\dd  t =
\int_{-1}^{\min\{\langle \bsx,\bsz\rangle, \langle \bsy,\bsz
\rangle\}} \,\dd  t = 1+ \min\{\langle \bsx , \bsz \rangle,
\langle \bsy, \bsz\rangle \}.
\end{equation*}
Thus
\begin{align*}
K_{\mathscr{C}}(\bsx,\bsy)
&= \int_{-1}^1 \int_{\mathbb{S}^d} 1_{C(\bsz;t)}(\bsx) 1_{C(\bsz;t)}(\bsy) \,\dd \sigma_d(\bsz) \,\dd  t \\
&= 1 + \int_{\mathbb{S}^d} \min\{ \langle \bsx, \bsz\rangle , \langle \bsy,\bsz \rangle\} \,\dd  \sigma_d(\bsz).
\end{align*}
We have 
\begin{equation*}
\min\{ \langle \bsx, \bsz\rangle , \langle \bsy,\bsz \rangle\} = \frac{1}{2} \left[\langle \bsx,\bsz \rangle + \langle \bsy, \bsz
\rangle - \left|\langle \bsx-\bsy, \bsz \rangle \right|\right]
\end{equation*}
and
\begin{equation*}
\int_{\mathbb{S}^d} \langle \bsx, \bsz \rangle \,\dd  \sigma_d(\bsz) = 0.
\end{equation*}
If $\bsx=\bsy$, then we therefore obtain 
\begin{equation*}
\int_{\mathbb{S}^d} \min\{ \langle \bsx, \bsz \rangle , \langle \bsy, \bsz \rangle\} \,\dd  \sigma_d(\bsz) = 0.
\end{equation*}
Let now $\bsx \neq \bsy$. Then
\begin{align*}
\int_{\mathbb{S}^d} \min\{ \langle \bsx, \bsz \rangle , \langle \bsy, \bsz \rangle\} \,\dd  \sigma_d(\bsz) 
&= -\frac{1}{2} \int_{\mathbb{S}^d} \left|\langle \bsx-\bsy, \bsz \rangle \right| \,\dd  \sigma_d(\bsz) \\
&= -\left\|\bsx-\bsy \right\| \frac{1}{2} \int_{\mathbb{S}^d} \left| \left\langle \frac{\bsx-\bsy}{\|\bsx-\bsy\|}, \bsz \right \rangle \right|\,\dd  \sigma_d(\bsz).
\end{align*}
The last integral does not depend on the unit vector $(\bsx-\bsy)/\|\bsx-\bsy\|$ by rotational symmetry. Thus, we have (cf. Appendix~\ref{sec:aux.results})
\begin{equation} \label{eq:C.d}
C_d := \frac{1}{2} \int_{\mathbb{S}^d} \left| \langle \bsp, \bsz \rangle \right| \dd  \sigma_d(\bsz) = \frac{1}{d} \frac{\omega_{d-1}}{\omega_{d}} = \frac{\mathcal{H}_d(\mathbb{B}^d)}{\mathcal{H}_d(\mathbb{S}^d)} = \frac{1}{d} \frac{\Gamma((d+1)/2)}{\sqrt{\pi} \Gamma(d/2)} \sim \frac{1}{\sqrt{2 \pi \, d}} \quad \text{as $d \to \infty$.}
\end{equation}
(Here, $\omega_d$ is the surface area of $\mathbb{S}^d$, $\mathbb{B}^d$ is the unit ball in $\mathbb{R}^d$, $\mathcal{H}_d$ is the $d$-dimensional Hausdorff measure normalized such that the $d$-dimensional unit cube $[0,1)^d$ has measure one and $\Gamma(z)$ is the Gamma function.) 
%
%

Therefore we obtain the following closed form representation:
\begin{equation} \label{eq:closed.form}
K_{\mathscr{C}}(\bsx,\bsy) = 1 - C_d \left\|\bsx-\bsy \right\|.
\end{equation}
The reproducing kernel has the following properties: for $\bsy, \bsz \in \mathbb{S}^d$ we have
\begin{itemize}
\item $1-2C_d \leq K(\bsx,\bsy) \leq 1$;
\item $K(\bsx,\bsy) = 1-2C_d$ if and only if $\bsx=-\bsy$ (clearly, $0 < 2 C_d < 1$);
\item $K(\bsx,\bsy) = 1$ if and only if $\bsx=\bsy$;
\end{itemize}

Note that the Karhunen-Loevy expansion of the function
$\|\bsx-\bsy\|$ is based on ultraspherical harmonics. Hence the
eigenfunctions of $K_{\mathscr{C}}$ are the ultraspherical
harmonics. The corresponding eigenvalues are also known. Therefore,
the functions in $\mathcal{H}_{\mathscr{C}}$ can be expanded using
ultraspherical harmonics and the inner product can be written using
the coefficients of such an expansion. See \cite{BrWo2010d_pre} for these
results.

\section{Worst-case error}

Let $\|f\|_{K_{\mathscr{C}}} = \sqrt{\langle f, f \rangle_{K_{\mathscr{C}}}}$ denote the norm in $\mathcal{H}_{\mathscr{C}}$. Then we define the worst-case error for a quadrature rule $Q_N$ given in \eqref{eq:Q.N} by
\begin{equation*}
e(\mathcal{H}_{\mathscr{C}}, Q_N) = \sup \left\{ \left|
\int_{\mathbb{S}^d} f(\bsx) \,\dd  \bsx - Q_N(f) \right| : f \in \mathcal{H}_{\mathscr{C}}, \|f\|_{K_{\mathscr{C}}} \le 1 \right\}.
\end{equation*}
Let $f \in \mathcal{H}_{\mathscr{C}}$. Then, by the reproducing kernel property $f(\bsy) = \langle f, K_{\mathscr{C}}(\cdot, \bsy) \rangle_{K_{\mathscr{C}}}$ for $\bsy \in \mathbb{S}^d$, 
and, since the integration functional $f \mapsto \int_{\mathbb{S}^d} f \dd \sigma_d$ is bounded on $\mathcal{H}_{\mathscr{C}}$ and has the 'representer' $\int_{\mathbb{S}^d} K_{\mathscr{C}}(\cdot, \bsz) \dd \sigma_d(\bsz)$, one can write
\begin{equation*}
\int_{\mathbb{S}^d} f(\bsx) \,\dd  \sigma_d(\bsx) - \frac{1}{N} \sum_{k=0}^{N-1} f(\bsx_k) = \left\langle f, \mathcal{R}(\mathcal{H}_{\mathscr{C}}, Q_N; \cdot) \right\rangle_{K_{\mathscr{C}}},
\end{equation*}
where the 'representer' of the error of numerical integration for the rule $Q_N$ for functions in $\mathcal{H}_{\mathscr{C}}$ is given by
\begin{equation*}
\mathcal{R}(\mathcal{H}_{\mathscr{C}}, Q_N; \bsx) = \int_{\mathbb{S}^d} K_{\mathscr{C}}(\bsx, \bsy) \,\dd  \sigma_d(\bsy) - \frac{1}{N} \sum_{k=0}^{N-1} K_{\mathscr{C}}(\bsx, \bsx_k), \qquad \bsx \in \mathbb{S}^d.
\end{equation*}
The Cauchy-Schwarz inequality yields
\begin{equation*}
\left| \int_{\mathbb{S}^d} f(\bsx) \,\dd  \sigma_d(\bsx) - \frac{1}{N} \sum_{k=0}^{N-1} f(\bsx_k) \right| = \left| \left\langle f, \mathcal{R}(\mathcal{H}_{\mathscr{C}}, Q_N; \cdot) \right\rangle_{K_{\mathscr{C}}} \right| \leq \|f\|_{K_{\mathscr{C}}} \left\| \mathcal{R}(\mathcal{H}_{\mathscr{C}}, Q_N; \cdot) \right\|_{K_{\mathscr{C}}}.
\end{equation*}
In particular, equality is assumed in the last relation when taking $f$ to be the 'representer' $\mathcal{R}(\mathcal{H}_{\mathscr{C}}, Q_N; \cdot)$ itself. It follows that
\begin{equation} \label{eq:e}
e(\mathcal{H}_{\mathscr{C}}, Q_N) = \left\| \mathcal{R}(\mathcal{H}_{\mathscr{C}}, Q_N; \cdot) \right\|_{K_{\mathscr{C}}} = \left\|\int_{\mathbb{S}^d} K_{\mathscr{C}}(\cdot, \bsx) \,\dd  \sigma_d(\bsx) - \frac{1}{N}
\sum_{k=0}^{N-1} K_{\mathscr{C}}(\cdot, \bsx_k) \right\|_{K_{\mathscr{C}}}.
\end{equation}
Expanding the square of the worst-case error and substituting the closed form of the reproducing kernel we arrive at the well-known representation
\begin{align}
&\left[ e(\mathcal{H}_{\mathscr{C}}, Q_N) \right]^2 = \left\langle \mathcal{R}(\mathcal{H}_{\mathscr{C}}, Q_N; \cdot), \mathcal{R}(\mathcal{H}_{\mathscr{C}}, Q_N; \cdot) \right\rangle \notag \\
&\phantom{eq}= \int_{\mathbb{S}^d} \int_{\mathbb{S}^d} K_{\mathscr{C}}(\bsx,\bsy) \,\dd \sigma_d(\bsx) \,\dd \sigma_d(\bsy) - \frac{2}{N} \sum_{k=0}^{N-1} \int_{\mathbb{S}^d} K_{\mathscr{C}}(\bsx,\bsx_k) \,\dd \sigma_d(\bsx) + \frac{1}{N^2} \sum_{k,\ell=0}^{N-1} K_{\mathscr{C}}(\bsx_k,\bsx_\ell) \notag \\
&\phantom{eq}= C_d \left[ \int_{\mathbb{S}^d} \int_{\mathbb{S}^d} \|\bsx-\bsy\|\,\dd \sigma_d(\bsx) \,\dd \sigma_d(\bsy) - \frac{1}{N^2} \sum_{k,\ell=0}^{N-1} \|\bsx_k - \bsx_\ell \| \right]. \label{eq:stolarsky.id.half1}
\end{align}

This shows how the square worst-case error in our reproducing kernel Hilbert space is related to \eqref{eq_intdoublesum}. In the next section we show how the worst-case error  $e(\mathcal{H}_{\mathscr{C}}, Q_N)$ in our reproducing kernel Hilbert space relates to the $L_2$ spherical cap discrepancy.

\section{Spherical cap discrepancy and Stolarsky's invariance principle}

Using the integral representation of the reproducing kernel \eqref{eq:reproducing.kernel} we have
\begin{align*}
\int_{\mathbb{S}^d} K_{\mathscr{C}}(\bsx,\bsy) \, \dd \sigma_d(\bsy) &= \int_{-1}^1 \int_{\mathbb{S}^d} 1_{C(\bsz;t)}(\bsx) \sigma_d({C(\bsz;t)}) \,\dd \sigma_d(\bsz) \,\dd t
\intertext{and}
\frac{1}{N} \sum_{k=0}^{N-1} K_{\mathscr{C}}(\bsx,\bsx_k) &= \int_{-1}^1 \int_{\mathbb{S}^d} \frac{1}{N} \sum_{k=0}^{N-1} 1_{C(\bsx;t)}(\bsz) 1_{C(\bsx_k;t)}(\bsz) \,\dd \sigma_d(\bsz) \,\dd t.
\end{align*}
Thus, the 'representer' of the error of numerical integration is of the form \eqref{f_rep}; that is
\begin{equation*}
\mathcal{R}(\mathcal{H}_{\mathscr{C}}, Q_N; \bsx) 
= \int_{-1}^1 \int_{\mathbb{S}^d} 1_{C(\bsz;t)}(\bsx) \left[ \sigma_d({C(\bsz;t)}) - \frac{1}{N} \sum_{k=0}^{N-1} 1_{C(\bsx_k;t)}(\bsz) \right] \dd \sigma_d(\bsz) \,\dd t.
\end{equation*}
Therefore, using the inner product representation \eqref{inner_product} in \eqref{eq:e}, we obtain
\begin{equation} \label{eq:ell.two.discr.id}
e(\mathcal{H}_{\mathscr{C}}, Q_N) = \left(\int_{-1}^1 \int_{\mathbb{S}^d} \left|\sigma_d(C(\bsz;t)) - \frac{1}{N} \sum_{k=0}^{N-1} 1_{C(\bsz;t)}(\bsx_k) \right|^2 \dd  \sigma_d(\bsz) \,\dd  t \right)^{1/2}
\end{equation}
stating that the worst-case error of the numerical integration formula $Q_N$ in \eqref{eq:Q.N} in the considered Sobolev space setting equals the so-called spherical cap $\mathbb{L}_2$-discrepancy of the integration nodes.

Combining \eqref{eq:stolarsky.id.half1} and \eqref{eq:ell.two.discr.id}, we arrive at Stolarsky's invariance principle for the Euclidean distance on spheres.

\begin{proposition}[Stolarsky~\cite{St1973}]
Let $\bsx_0,\ldots, \bsx_{N-1} \in \mathbb{S}^d$ be an arbitrary $N$ point configuration on the sphere $\mathbb{S}^d$. Then we have
\begin{equation}
\begin{split} \label{eq:Stolarsky.inv.principle}
&\frac{1}{N^2} \sum_{k,\ell=0}^{N-1} \|\bsx_k - \bsx_\ell \| + \frac{1}{C_d} \int_{-1}^1 \int_{\mathbb{S}^d} \left|\sigma_d(C(\bsz;t)) - \frac{1}{N} \sum_{k=0}^{N-1} 1_{C(\bsz;t)}(\bsx_k) \right|^2 \dd  \sigma_d(\bsz) \,\dd  t \\
&\phantom{equals}= \int_{\mathbb{S}^d} \int_{\mathbb{S}^d} \|\bsx-\bsy\|\,\dd \sigma_d(\bsx) \,\dd \sigma_d(\bsy).
\end{split}
\end{equation}
\end{proposition}

The $\mathbb{L}_2$-discrepancy of an $N$-point configuration on $\mathbb{S}^d$ decreases as its sum of distances increases and vice versa. The right-hand side is the distance integral of the uniform measure $\sigma_d$ on the sphere $\mathbb{S}^d$ which is the unique extremal measure (also known as the equilibrium measure) maximizing the distance integral
\begin{equation*}
\mathcal{I}[\mu] := \int_{\mathbb{S}^d} \int_{\mathbb{S}^d} \left\| \bsx - \bsy \right\| \dd \mu( \bsx ) \dd \mu( \bsy )
\end{equation*}
over the family of (Borel) probability measures $\mu$ supported on $\mathbb{S}^d$. For the potential theory of the generalized distance integral we refer to Bj{\"o}rck \cite{Bj1956}.

\section{A weighted reproducing kernel}

The above results can be generalized by introducing a weight
function. Let $v:[-1,1] \to\mathbb{R}$ satisfy $v(t) > 0$ for all
$t$ and which has an antiderivative, which we denote by $V$. Then we define the
reproducing kernel with weight function $v$ as follows
\begin{equation}\label{weighted_kernel}
K_{\mathscr{C},v}(\bsx,\bsy) = \int_{-1}^1 v(t) \int_{\mathbb{S}^d} 1_{C(\bsz;t)}(\bsx) 1_{C(\bsz;t)}(\bsy) \,\dd  \sigma_d(\bsz) \,\dd  t, \qquad \bsx, \bsy \in \mathbb{S}^d.
\end{equation}
For functions represented by integrals
\begin{equation*} 
f_i(\bsx) = \int_{-1}^1 \int_{\mathbb{S}^d} g_i(\bsz;t) 1_{C(\bsz;t)}(\bsx) \,\dd \sigma_d(\bsz) \,\dd  t, \qquad i = 1, 2,
\end{equation*}
the corresponding inner product is now given by
\begin{equation*}
\langle f_1, f_2 \rangle_{K_{\mathscr{C}, v}} = \int_{-1}^1 \frac{1}{v(t)} \int_{\mathbb{S}^d} g_1(\bsz;t) \overline{g_2(\bsz;t)} \,\dd  \sigma_d(\bsz) \,\dd  t.
\end{equation*}

The reproducing kernel can be written as
\begin{equation} \label{eq:weighted.repr.kernel}
K_{\mathscr{C}, v}(\bsx,\bsy) = \int_{\mathbb{S}^d} V(\min\{\langle \bsx, \bsz \rangle, \langle \bsy, \bsz \rangle\}) \,\dd  \sigma_d(\bsz) - V(-1), \qquad \bsx, \bsy \in \mathbb{S}^d.
\end{equation}
For certain weight functions $v$,  this expression may have a
concise form. This reproducing kernel defines a reproducing kernel
Hilbert space $\mathcal{H}_{\mathscr{C},v}$.

The 'representer' of the error of numerical integration for the rule $Q_N$ for functions in $\mathcal{H}_{\mathscr{C},v}$ takes on the form
\begin{equation*}
\mathcal{R}(\mathcal{H}_{\mathscr{C}, v}, Q_N; \bsx) = \int_{-1}^1 \int_{\mathbb{S}^d} 1_{C(\bsz;t)}(\bsx) \, v(t) \left[ \sigma_d({C(\bsz;t)}) - \frac{1}{N} \sum_{k=0}^{N-1} 1_{C(\bsx_k;t)}(\bsz) \right] \dd \sigma_d(\bsz) \,\dd t.
\end{equation*}
We claim that $K_{\mathscr{C}, v}(\bsx,\bsy)$ is a function of the inner product $\langle \bsx, \bsy \rangle$, cf. Appendix~\ref{sec:app.weighted.repr.kernel}. Using the same approach as before we obtain
\begin{equation} \label{eq_ev1}
\left[ e(\mathcal{H}_{\mathscr{C},v}, Q_N) \right]^2 = \frac{1}{N^2} \sum_{k,\ell=0}^{N-1} K_{\mathscr{C},v}(\bsx_k,\bsx_\ell) - \int_{\mathbb{S}^d} \int_{\mathbb{S}^d} K_{\mathscr{C},v}(\bsx,\bsy) \,\dd \sigma_d(\bsx) \,\dd  \sigma_d(\bsy). 
\end{equation}
This worst case error can also be expressed in terms of a weighted discrepancy measure:
\begin{align}
e(\mathcal{H}_{\mathscr{C},v}, Q_N) 
&= \left\|\int_{\mathbb{S}^d} K_{\mathscr{C},v}(\cdot, \bsy) \,\dd  \sigma_d(\bsy) - \frac{1}{N} \sum_{k=0}^{N-1} K_{\mathscr{C},v}(\cdot, \bsx_k) \right\|_{K_{\mathscr{C},v}} = \left\| \mathcal{R}(\mathcal{H}_{\mathscr{C}, v}, Q_N; \cdot) \right\|_{K_{\mathscr{C},v}} \notag \\
&= \left(\int_{-1}^1 v(t) \int_{\mathbb{S}^d} \left|\sigma_d(C(\bsx;t)) - \frac{1}{N} \sum_{k=0}^{N-1} 1_{C(\bsx;t)}(\bsx_k) \right|^2 \,\dd  \sigma_d(\bsx) \,\dd  t \right)^{1/2}. \label{eq_ev2}
\end{align}
Using \eqref{eq_ev1} and \eqref{eq_ev2} we obtain the weighted version of Stolarsky invariance principle.

\begin{theorem} \label{thm:weighted.kernel}
Let $\bsx_0,\ldots, \bsx_{N-1} \in \mathbb{S}^d$ be an arbitrary $N$ point configuration on the sphere $\mathbb{S}^d$. Let $K_{\mathscr{C},v}$ be the weighted reproducing kernel given by \eqref{weighted_kernel}. Then we have
\begin{equation*}
\begin{split}
&\frac{1}{N^2} \sum_{k,\ell=0}^{N-1} K_{\mathscr{C},v}(\bsx_k,\bsx_\ell) + \int_{-1}^1  v(t) \int_{\mathbb{S}^d} \left| \sigma_d(C(\bsx;t)) - \frac{1}{N} \sum_{k=0}^{N-1} 1_{C(\bsx;t)}(\bsx_k) \right|^2 \,\mathrm{d} \sigma_d(\bsx) \,\mathrm{d} t \\   &\phantom{equals}= \int_{\mathbb{S}^d} \int_{\mathbb{S}^d} K_{\mathscr{C},v}(\bsx,\bsy) \,\mathrm{d} \sigma_d(\bsx) \,\mathrm{d} \sigma_d(\bsy).
\end{split}
\end{equation*}
\end{theorem}

The double integral above can be expressed in terms of the weight function, see \eqref{eq:double.int}.

In \cite{Br2003}, Stolarsky's (general) invariance principle is extended and used to get bounds for the spherical cap discrepancy, see also \cite{Br2004}. Stolarksy~\cite{St1975b} also extended his principle to certain metric spaces arising from measures.

Stolarsky~\cite{St1973} introduced the function
\begin{equation*}
\rho(\bsx, \bsy) = \int_{\mathbb{S}^d} \int_{\min\{ \langle \bsx, \bsy \rangle, \langle \bsy, \bst \rangle \}}^{\max\{ \langle \bsx, \bsy \rangle, \langle \bsy, \bst \rangle \}} g(u) \dd u \dd \sigma_d( \bsz ), \qquad \bsx, \bsy \in \mathbb{S}^d,
\end{equation*}
which becomes a metric if the kernel $g$ (integrable on $[0,1]$) is positive but the proof of the corresponding invariance principle
\begin{equation} 
\begin{split} \label{eq:Stolarksy.general.principle}
&\frac{1}{N^2} \sum_{k,\ell=0}^{N-1} \rho(\bsx_k,\bsx_\ell) + 2 \int_{-1}^1 v(t) \int_{\mathbb{S}^d} \left|\sigma_d(C(\bsx;t)) - \frac{1}{N} \sum_{k=0}^{N-1} 1_{C(\bsx;t)}(\bsx_k) \right|^2 \,\dd  \sigma_d(\bsx) \,\dd  t \\
&\phantom{equals}= \int_{\mathbb{S}^d} \int_{\mathbb{S}^d} \rho(\bsx,\bsy) \,\dd \sigma_d(\bsx) \,\dd  \sigma_d(\bsy)
\end{split}
\end{equation}
making use of Haar integrals over the special orthogonal group $\mathrm{SO}(d+1)$ does not require it. Note that for $g \equiv 1$ the function $\rho( \bsx, \bsy )$ is a constant multiple of the Euclidean distance. With some care one may even consider $g(x) = 1 / ( 1 - x^2 )$.

It is well-known that a reproducing kernel $K(\bsx,\bsy)$ induces a distance (metric) by means of
\begin{equation*}
\dist( \bsx, \bsy ) = \dist_{K}( \bsx, \bsy ) = \sqrt{K(\bsx,\bsx) - 2 K(\bsx,\bsy) + K(\bsy,\bsy)}. 
\end{equation*}
For example, the reproducing kernel \eqref{eq:closed.form} yields
\begin{equation*}
\dist_{K_{\mathscr{C}}}( \bsx, \bsy ) = \sqrt{C_d} \, \sqrt{ \left\| \bsx - \bsy \right\| }.
\end{equation*}
In general, for the symmetric weighted kernel $K_{\mathscr{C},v}(\bsx, \bsy)$, which does only depend on the inner product $\langle \bsx, \bsy \rangle$, it follows that ($\bsa \in \mathbb{S}^d$ fixed)
\begin{equation*}
K_{\mathscr{C},v}(\bsx, \bsy) = \frac{1}{2} \left\{ \left[ \dist_{\mathscr{C},v}( \bsx, \bsy ) \right]^2 - K_{\mathscr{C},v}(\bsx, \bsx) - K_{\mathscr{C},v}(\bsy, \bsy) \right\} = \frac{1}{2} \left[ \dist_{\mathscr{C},v}( \bsx, \bsy ) \right]^2 - K_{\mathscr{C},v}(\bsa, \bsa).
\end{equation*}
By Theorem~\ref{thm:weighted.kernel} on arrives at
\begin{equation*}
\begin{split}
&\frac{1}{N^2} \sum_{k,\ell=0}^{N-1} \left[ \dist_{\mathscr{C},v}( \bsx, \bsy ) \right]^2 + 2 \int_{-1}^1  v(t) \int_{\mathbb{S}^d} \left| \sigma_d(C(\bsx;t)) - \frac{1}{N} \sum_{k=0}^{N-1} 1_{C(\bsx;t)}(\bsx_k) \right|^2 \,\mathrm{d} \sigma_d(\bsx) \,\mathrm{d} t \\   
&\phantom{equals}= \int_{\mathbb{S}^d} \int_{\mathbb{S}^d} \left[ \dist_{\mathscr{C},v}( \bsx, \bsy ) \right]^2 \,\mathrm{d} \sigma_d(\bsx) \,\mathrm{d} \sigma_d(\bsy),
\end{split}
\end{equation*}
which should be compared with \eqref{eq:Stolarksy.general.principle}.

\vspace{10mm}
{\bf Acknowledgement:} The first author is grateful to the School of Mathematics and Statistics at UNSW for their support.

\appendix

\section{Auxiliary results}
\label{sec:aux.results}

The normalized surface area measure $\sigma_d$ on $\mathbb{S}^d$ admits the following decomposition
\begin{equation} \label{eq:sigma.d}
\dd \sigma_d( \bsy ) = \frac{\omega_{d-1}}{\omega_d} \left( 1 - t^2 \right)^{d/2-1} \dd t \dd \sigma_{d-1}( \bsy^*), \qquad \bsy = ( \sqrt{1-t^2} \, \bsy^*, t) \in \mathbb{S}^d,
\end{equation}
where $t \in [-1,1]$, $\bsy^* \in \mathbb{S}^{d-1}$ and $\omega_d$ denotes the surface area of $\mathbb{S}^d$ (cf M{\"u}ller~\cite{Mu1966}). (By definition $\langle \bsy, \bsp \rangle = t$, where $\bsp$ is the North Pole of $\mathbb{S}^d$.) Thus, by rotational symmetry, the integral of a zonal function $f(\langle \bsz, \boldsymbol{\cdot} \rangle)$, $\bsz \in \mathbb{S}^d$ fixed, with respect to $\sigma_d$ reduces to 
\begin{equation*}
\int_{\mathbb{S}^d} f(\langle \bsz, \bsy \rangle) \dd \sigma_d( \bsy ) = \int_{\mathbb{S}^d} f(\langle \bsp, \bsy \rangle) \dd \sigma_d( \bsy ) = \frac{\omega_{d-1}}{\omega_d} \int_{-1}^1 f(t) \left( 1 - t^2 \right)^{d/2-1} \dd t.
\end{equation*}

\begin{proof}[Proof of relations \eqref{eq:C.d}]
One gets
\begin{align*}
C_d 
&= \frac{1}{2} \int_{\mathbb{S}^d} \left| \langle \bsp, \bsy \rangle \right| \dd \sigma_d(\bsy) = \frac{1}{2} \frac{\omega_{d-1}}{\omega_d} \int_{-1}^1 \left| t \right| \left( 1 - t^2 \right)^{d/2-1} \dd t = \frac{1}{2} \frac{\omega_{d-1}}{\omega_d}  \int_{0}^1 \left( 1 - t^2 \right)^{d/2-1} 2 t \dd t \\
&= \frac{1}{d} \frac{\omega_{d-1}}{\omega_d} = \frac{1}{d} \frac{\Gamma((d+1)/2)}{\sqrt{\pi} \Gamma(d/2)} = \frac{\mathcal{H}_d(\mathbb{B}^d)}{\mathcal{H}_d(\mathbb{S}^d)} \sim \frac{1}{\sqrt{2 \pi \, d}} \quad \text{as $d \to \infty$.}
\end{align*}
The second equality follows from 
\begin{equation*}
1 = \sigma_d( \mathbb{S}^d ) = \frac{\omega_{d-1}}{\omega_d} \int_{-1}^1 \left( 1 - t^2 \right)^{d/2-1} \dd t = \frac{\omega_{d-1}}{\omega_d} \int_0^1 v^{1/2-1} \left( 1 - v \right)^{d/2-1}  \dd v = \frac{\omega_{d-1}}{\omega_d} \betafcn(1/2,d/2),
\end{equation*}
where $\betafcn(a,b) = \Gamma(a) \Gamma(b) / \Gamma(a+b)$ is the beta function. The third equality follows from the well-known formulas for the volume of the unit ball in $\mathbb{R}^d$ and the surface area of $\mathbb{S}^d$. The asymptotics follows from the asymptotic expansion of a ratio of Gamma functions (cf. \cite{DLMF2010.05.07}).
\end{proof}

\section{The weighted reproducing kernel}
\label{sec:app.weighted.repr.kernel}

Next, we investigate the weighted reproducing kernel \eqref{eq:weighted.repr.kernel} in more detail. In particular, it will be shown that the kernel $K_{\mathscr{C}, v}(\bsx,\bsy)$ is a function of the inner product $\langle \bsx, \bsy \rangle$.

On observing that $\langle \bsx, \bsz \rangle \leq \langle \bsy, \bsz \rangle$ if and only if $\langle \bsy - \bsx, \bsz \rangle \geq 0$ we may write for $\bsx \neq \bsy$ 
\begin{equation*}
K_{\mathscr{C}, v}(\bsx,\bsy) = \mathcal{A}_{\mathscr{C}, v}(\bsx,\bsy) + \mathcal{A}_{\mathscr{C}, v}(\bsy,\bsx) - V(-1), 
\end{equation*}
which immediately shows symmetry of the reproducing kernel, where 
\begin{equation} \label{eq:A.C.v}
\mathcal{A}_{\mathscr{C}, v}(\bsx,\bsy) = \int_{\mathbb{S}^d} V( \langle \bsx, \bsz \rangle ) \, 1_{[0,1]}( \left\langle \frac{\bsy - \bsx}{\| \bsy - \bsx \|}, \bsz \right\rangle ) \dd \sigma_d( \bsz ).
\end{equation}
By abuse of notation we set (note that $\langle \bsx, \bsy \rangle = u$)
\begin{align*}
\bsz &= t \, \bsx + \sqrt{1-t^2} \, \bsz^*, \qquad -1 \leq t \leq 1, \bsz^* \in \mathbb{S}^{d-1}, \\
\bsy &= u \, \bsx + \sqrt{1-u^2} \, \bsy^*, \qquad -1 \leq t \leq 1, \bsy^* \in \mathbb{S}^{d-1}.
\end{align*}
In this way $\bsx$ will be the 'North Pole' in the decomposition \eqref{eq:sigma.d} and we obtain
\begin{equation*}
\mathcal{A}_{\mathscr{C}, v}(\bsx,\bsy) = \frac{\omega_{d-1}}{\omega_d} \int_{-1}^1 V(t) \left\{ \int_{\mathbb{S}^{d-1}} 1_{[0,1]}(\left\langle \frac{\bsy - \bsx}{\| \bsy - \bsx \|}, \bsz \right\rangle) \sigma_{d-1}( \bsz^* ) \right\} \left( 1 - t^2 \right)^{d/2-1} \dd t.
\end{equation*}
The indicator function in the inner integral is a zonal function depending on $w = \langle \bsy^*, \bsz^* \rangle$ only. Thus, we apply again \eqref{eq:sigma.d} with $\bsy^*$ as 'North Pole'. That is
\begin{equation*}
\mathcal{A}_{\mathscr{C}, v}(\bsx,\bsy) = \frac{\omega_{d-1}}{\omega_d} \int_{-1}^1 V(t) \left\{ \frac{\omega_{d-2}}{\omega_{d-1}} \int_{-1}^1 1_{[0,1]}(\left\langle \frac{\bsy - \bsx}{\| \bsy - \bsx \|}, \bsz \right\rangle) \left( 1 - w^2 \right)^{(d-1)/2-1} \dd w \right\} \left( 1 - t^2 \right)^{d/2-1} \dd t, 
\end{equation*}
where the inner product evaluates as
\begin{equation} \label{eq:arg.indicator.fcn}
\left\langle \frac{\bsy - \bsx}{\| \bsy - \bsx \|}, \bsz \right\rangle = \sqrt{1-t^2} \, \sqrt{\frac{1+\langle \bsx, \bsy \rangle}{2}} \, w - \sqrt{\frac{1-\langle \bsx, \bsy \rangle}{2}} \, t, \qquad \bsx, \bsy, \bsz \in \mathbb{S}^d, \bsx \neq \bsy.
\end{equation}
Proceeding similarly for $\mathcal{A}_{\mathscr{C}, v}(\bsy,\bsx)$, one sees that, indeed, $\mathcal{A}_{\mathscr{C}, v}(\bsy,\bsx) = \mathcal{A}_{\mathscr{C}, v}(\bsx,\bsy)$. Furthermore, $\mathcal{A}_{\mathscr{C}, v}(\bsx,\bsy)$ depends only on the inner product $\langle \bsx, \bsy \rangle$ which in turn implies that the reproducing kernel $K_{\mathscr{C}, v}(\bsx,\bsy)$ is a function of the inner product $\langle \bsx, \bsy \rangle$.

The right-hand side in \eqref{eq:arg.indicator.fcn} describes a line which stays strictly between the levels $-1$ and $1$ for $v$ in $[-1,1]$ by the left-hand side in \eqref{eq:arg.indicator.fcn}. Further analysis gives that the indicator functions is one (i) if $t \leq - \sqrt{(1+u)/2}$ and $-1\leq v \leq 1$, or, (ii) if $- \sqrt{(1+u)/2} \leq t \leq \sqrt{(1+u)/2}$ and $( t / \sqrt{1-t^2} ) \sqrt{(1-u)/(1+u)} \leq v \leq 1$, and zero otherwise. This leads to
\begin{equation*}
\begin{split}
\mathcal{A}_{\mathscr{C}, v}(\bsx,\bsy) 
&= \frac{\omega_{d-1}}{\omega_d} \int_{-1}^{-\sqrt{(1+u)/2}} V(t) \left( 1 - t^2 \right)^{d/2-1} \dd t \\
&\phantom{=}+ \frac{\omega_{d-1}}{\omega_d} \int_{-\sqrt{(1+u)/2}}^{\sqrt{(1+u)/2}} V(t) \betafcnReg_{(1-x(t))/2}((d-1)/2,(d-1)/2) \left( 1 - t^2 \right)^{d/2-1} \dd t,
\end{split}
\end{equation*}
where, when using $u = \langle \bsx, \bsy \rangle = \cos \phi$ ($0 < \phi < \pi$) and $t = \cos \psi$, one has
\begin{equation*}
\sqrt{\frac{1+u}{2}} = \cos(\phi/2), \qquad \sqrt{\frac{1-u}{1+u}} = \tan(\phi/2), \qquad x(t) = \sqrt{\frac{1-u}{1+u}} \, \frac{t}{\sqrt{1-t^2}} = \frac{\cot \psi}{\cot( \phi/2)}.
\end{equation*}
The change of variable $\xi = x(t)$ yields
\begin{equation*}
\begin{split}
\mathcal{A}_{\mathscr{C}, v}(\bsx,\bsy) 
&= \frac{\omega_{d-1}}{\omega_d} \int_{-1}^{-\sqrt{(1+u)/2}} V(t) \left( 1 - t^2 \right)^{d/2-1} \dd t \\
&\phantom{=}+ \frac{\omega_{d-1}}{\omega_d} \left( \frac{1-u}{1+u} \right)^{d/2} \int_{-1}^{1} V(\frac{\xi}{\sqrt{\frac{1-u}{1+u}+\xi^2}}) \betafcnReg_{(1-\xi)/2}((d-1)/2,(d-1)/2) \frac{\dd \xi}{\left( \frac{1-u}{1+u}+\xi^2 \right)^{(d+1)/2}},
\end{split}
\end{equation*}
where we make use of the regularized incomplete beta function
\begin{equation*}
\betafcnReg_z(a,b) = \betafcn_z(a,b) / \betafcn(a,b), \qquad \betafcn_z(a,b) = \int_0^z t^{a-1} (1 - t)^{b-1} \dd t, \quad a,b > 0.
\end{equation*}

We compute the following integral (using \eqref{eq:A.C.v}):
\begin{equation*}
\int_{\mathbb{S}^d} \mathcal{A}_{\mathscr{C}, v}(\bsx,\bsy) \, \dd \sigma_d( \bsy ) = \int_{\mathbb{S}^d} V( \langle \bsx, \bsz \rangle ) \, \int_{\mathbb{S}^d} 1_{[0,1]}( \left\langle \frac{\bsy - \bsx}{\| \bsy - \bsx \|}, \bsz \right\rangle )  \, \dd \sigma_d( \bsy ) \, \dd \sigma_d( \bsz ).
\end{equation*}
The inner integral is one if $\bsz$ is in the half-sphere centered at $-\bsx$ and zero otherwise. Hence,
\begin{equation*}
\int_{\mathbb{S}^d} \mathcal{A}_{\mathscr{C}, v}(\bsx,\bsy) \, \dd \sigma_d( \bsy ) = \frac{\omega_{d-1}}{\omega_d} \int_0^1 V(-t) \left( 1 - t^2 \right)^{d/2-1} \dd t.
\end{equation*}
Since $g^\prime(t) = ( \omega_{d-1} / \omega_d ) (1 - t^2)^{d/2-1}$ for $g(t) = (1/2) \betafcnReg_{t^2}(1/2,d/2)$, $0 \leq t \leq 1$, integration by parts gives
\begin{equation*}
\int_{\mathbb{S}^d} \mathcal{A}_{\mathscr{C}, v}(\bsx,\bsy) \, \dd \sigma_d( \bsy ) = \frac{1}{2} V(-1) + \frac{1}{2} \int_0^1 v(-t) \betafcnReg_{t^2}(1/2,d/2) \dd t.
\end{equation*}
It follows that
\begin{subequations} \label{eq:double.int}
\begin{align}
\int_{\mathbb{S}^d} \int_{\mathbb{S}^d} K_{\mathscr{C}, v}(\bsx,\bsy) \, \dd \sigma_d( \bsx ) \, \dd \sigma_d( \bsy ) 
&= 2 \frac{\omega_{d-1}}{\omega_d} \int_0^1 V(-t) \left( 1 - t^2 \right)^{d/2-1} \dd t - V(-1) \\
&= \int_0^1 v(-t) \betafcnReg_{t^2}(1/2,d/2) \dd t.
\end{align}
\end{subequations}

\bibliographystyle{abbrv}
\bibliography{bibliography}

\end{document}